\documentclass[a4paper,11pt]{article}
\usepackage{amssymb,amsmath}
\author{Florent Benaych-Georges} 
\title{Failure of the Raikov Theorem for Free Random Variables}
\newcommand{\D}{\mathrm{d}}
\newtheorem{theorem}{Theorem}
\newtheorem{lemma}[theorem]{Lemma}
\newtheorem{corollary}[theorem]{Corollary}
\newenvironment{proof}{\noindent {\bf Proof }}{\ \ \ $\square$}
\begin{document}
\maketitle
\begin{abstract} We show that the sum of two free random variables can have a free Poisson law without any of them having a free Poisson law.
\end{abstract}

The classical convolution $*$ of probability measures on $\mathbb{R}$  has an analogue in the theory of free probability,  the free convolution $\boxplus$ (see \cite{defconv}), and the structures of semi-groups defined by classical and free convolutions on the set of probability measures on $\mathbb{R}$  present many analogies. For example, we can define infinitely divisible laws in the same way  for both convolutions, and there exists a natural isomorphism between the semi-group of classical infinitely divisible laws and the semi-group of free infinitely divisible laws (see \cite{appenice}). This isomorphism transforms a Gaussian law into the semi-circle law with the same mean and the same variance, a Poisson law into the Marchenko-Pastur law with the same mean (see \cite{hiai}). Moreover, there is a deep correspondance between  limit theorems  for weak convergence of  sums of independent random variables and limit theorems  for weak convergence of sums of free random variables (see \cite{appenice}). It then seems natural to see how far the analogy goes.

Two already established properties of free convolution show that the isomorphism between the semi-group of classical infinitely divisible laws and the semi-group of free infinitely divisible laws cannot be extended to the set of probability  measures. The first one is the following (see \cite{spei2}): for every probability  measure $\mu$ on the real line, there exists a family   $(\mu_t)_{t\geq 1}$ of probability  measures such that  $\mu_1=\mu$ and for all  $s,t\geq 1$, one has $\mu_{s+t}=\mu_s\boxplus\mu_t$.  There is no similar result when we replace $\boxplus$ by $*$ for  $\mu= \frac{1}{2}\left(\delta_0+\delta_1\right)$. The second one is the failure of the Cram\'er theorem  for free random variables: the sum of two free random variables can be distributed according to the semi-circle law without any of them having a semi-circle law$-$the classical Cram\'er theorem  states that if the sum of two independent random variables is Gaussian, then each of them is Gaussian. 

We will show in this article that the analogue of Raikov's theorem for free probability  is false: the sum of two free random variables can have a free Poisson law (i.e. a Marchenko-Pastur law) without any of them (even after translation) having a free Poisson law$-$the classical Raikov theorem  states that if the sum of two independent random variables has, after translation,  a Poisson law, then each of them has, after translation,  a Poisson law. The proof is very close to the proof given in \cite{supercv} of the failure of the Cram\'er theorem  for free random variables. 

In order to prove our result, we need to discuss the analytic method for calculating free convolutions. We will restrict ourselves to compactly supported probability  measures on $\mathbb{R}$. Let $\mathbb{C}^+$ (resp. $\mathbb{C}^-$) denote the upper (resp. lower) half plane. Given $\mu$ a compactly supported measure on $\mathbb{R}$, consider the function  \begin{eqnarray*}G_\mu\; :\; \mathbb{C}^+&\to&\mathbb{C}^-\\
z&\mapsto&\int \frac{\D \mu(t)}{z-t}.
\end{eqnarray*}
This function, called Cauchy transform  of $\mu$, is analytic, and \\
$\underset{|z|\to\infty}{\lim} zG_\mu(z)=1$. Hence there exists  a function $K_\mu$ that is meromorphic in a neighborhood in $\mathbb{C}^-$ of zero, with a single pole at zero such that   $G_\mu(K_\mu(z))=z$ for $z$ close to zero. One can write $K_\mu(z)=\frac{1}{z}+R_\mu(z)$, where $R_\mu$ is an analytic function in a neighborhood of zero in $\mathbb{C}^-$ (i.e. $V\cap \mathbb{C}^-$ with $V$ a neighborhood of zero in $\mathbb{C}$), called the $R$-transform of $\mu$. It was shown in   \cite{voic0} (see also in \cite{voic1}) that $R_\mu$ determines $\mu$ and that  $R_{\mu\boxplus\nu}=R_\mu+R_\nu$. 

The $R$-transform of the free Poisson law with mean $\lambda$, that we will denote by $\mu_\lambda$, is given by the formula \[R_{\mu_\lambda}(z)=\frac{\lambda}{1-z}.\]
 
\begin{theorem}Consider $0<\lambda<1$,  $r\in]0,\frac{\sqrt{\lambda}}{1-\lambda}[$. Let $(f_n)_{n\in\mathbb{N}}$ be a sequence of analytic functions on  $\left\{z\; ;\; \left|z- \frac{1}{1-\lambda}\right| >r\right\}$, such that 
\\
(i) for all $n$, for all $z$, $f_n(\overline{z})=\overline{f_n(z)}$, 
\\
(ii) for all $n$, $\underset{\infty}{\lim} z f_n'(z)=\underset{\infty}{\lim}f_n(z)=0$,
\\\
(iii) the sequence $(f_n)$ converges uniformly to $0$ when $n$ tends to infinity.
\\
Then when $n$ is large enough, the function $z\mapsto\frac{\lambda}{1-z}+f_n(z)$ is the  $R$-transform of a compactly supported probability measure on $\mathbb{R}$. 
\end{theorem}
Before giving the proof of the theorem,  here is the corollary that has led us to establish it.

 \begin{corollary}Consider $\lambda>0$. There exist two compactly supported probability measures on $\mathbb{R}$   who, even after translation, are not free Poisson laws, and whose   free convolution is the free Poisson law with mean $\lambda$.
\end{corollary}

\begin{proof}Consider $0<\lambda'<\min \{1, \lambda /\!2\}$.  It suffices to define, with  $\varepsilon$ small enough, $\mu_{\small{+}}$ to be the law whose  $R$-transform is $z\mapsto\frac{\lambda'}{1-z}+\frac{\varepsilon}{(z-1/\!(1-\lambda))^2}$ (such a law exists by the theorem), and  $\mu_{\small{-}}$ to be the law whose $R$-transform is $z\mapsto\frac{\lambda'}{1-z}-\frac{\varepsilon}{(z-1/\!(1-\lambda))^2}+\frac{\lambda-2\lambda'}{1-z}$ (such a law exists by the theorem  and because the sum of two $R$-transforms is a $R$-transform). Then the pair $\mu_{\small{+}}$, $\mu_{\small{-}}$ is suitable.
\end{proof}
As a preliminary to the proof of the theorem,  here is a lemma that states the existence of a measure  (see \cite{akhi}).
\begin{lemma}Let $G$ be an analytic function from the upper half plane 
$\mathbb{C}^+$
to the lower half plane 
$\mathbb{C}^-$ 
that satisfies:
\\
(i) $zG(z)$ tends to $1$ when $|z|$ tends to infinity,
\\
(ii) denoting by $K$ the inverse function of $G$ from a neighborhood of  $0$ in  $\mathbb{C}^-$ (i.e. from $V\cap \mathbb{C}^-$ with $V$ a neighborhood of zero in $\mathbb{C}$)  to a neighborhood of infinity in $\mathbb{C}^+$, the function $z\mapsto K(z)-\frac{1}{z}$ can be extended to a neighborhood of zero in $\mathbb{C}$.
\\
Then there exists a compactly supported probability measure on $\mathbb{R}$  whose Cauchy transform  is  $G$.
\end{lemma}
Here is the proof of the theorem.  $\Re$ and $\Im$ will denote respectively the real and imaginary parts.

\begin{proof} By the lemma, it is enough to show that when $n$ is large enough, there exists a domain $\Omega_n$ of $\mathbb{C}^-$, neighborhood of $0$ in $\mathbb{C}^-$,  such that  $K_n$ induces a bijection from $\Omega_n$ onto $\mathbb{C}^+$, where $K_n(z)=\frac{1}{z}+ \frac{\lambda}{1-z}+f_n(z)$; and that the inverse function of $K_n$ satisfies (i) and  (ii) of the lemma. For convenience, we will work with the  functions $\psi_n(z)=K_n\left(\frac{\sqrt{\lambda}}{1-\lambda}z+\frac{1}{1-\lambda}\right)$, and $\psi(z)=K_{\mu_\lambda}\left(\frac{\sqrt{\lambda}}{1-\lambda}z+\frac{1}{1-\lambda}\right)$.  Recall $K_{\mu_\lambda}(z)=\frac{1}{z}+ \frac{\lambda}{1-z}$.

Fix $\alpha \in ]0,1/\!10[$ such that  $1-2\alpha > \max\{r,\sqrt{\lambda}\}$ and $1+2\alpha <1/\!\sqrt{\lambda}$, and  fix $\eta\in ]0,\alpha/\!2[$. \\
First, let us note few properties of the function  $\psi$, that we will then extend to the functions  $\psi_n$ with $n$ large enough, using Rouché's theorem.  \\
 We have $$\psi(z)=\frac{z(1-\lambda)^2}{(\sqrt{\lambda}z+1)(z+\sqrt{\lambda})},$$
so zero is the unique preimage of zero by  $\psi$, and for every  non zero $Z$, the set of the preimages of $Z$ by $\psi$ is the set of the roots of  $$X^2+\left(\frac{1}{\sqrt{\lambda}}+\sqrt{\lambda}-\frac{(1-\lambda)^2}{Z\sqrt{\lambda}}\right)X+1.$$ Their product is one, so $\psi$ is one-to-one in $\mathbb{C}^+$ and in $\mathbb{C}^-$. We deduce:
(a) when $n$ is large enough, $\psi_n$ is one-to-one in a neighborhood of $A=\{ z\; ;\;  1-\eta\leq |z|\leq 1+\eta, |z\pm 1|\geq \eta, \Im z \leq 0\}$.
\\
For the same reason, $\psi$ is one-to-one in $\{z\in\mathbb{C}\; ;\; 0<|z|<1,\:z\neq-\sqrt{\lambda}\}$  and in $\{z\in\mathbb{C}\; ;\; 1<|z|,\:z\neq-\frac{1}{\sqrt{\lambda}}\}$. Hence: \\
(b) for $n$ large enough, $\psi_n$ is one-to-one in a neighborhood of the corona $B=\{z \; ;\;   1-2\eta\leq |z|\leq 1-\eta\}$ and in a neighborhood of the corona $C=\{z \; ;\;   1+\eta\leq |z|\leq 1+2\eta\}$.\\
We have $$\Im \left(\psi(z)\right)=\frac{(1-\lambda)^2\sqrt{\lambda}}{|\sqrt{\lambda}z+1||z+\sqrt{\lambda}|}(1-|z|^2)\Im z,$$
so $\psi(z)$ is real if and only if  $z$ is real or has modulus one, hence:\\
(c) for $n$ large enough,  $\Im \psi_n(z) <0$ for $z=re^{{\rm i}\theta}\in B$ with $\theta\in[-\pi+\eta,-\eta]$,\\
(d) for $n$ large enough,  $\Im \psi_n(z) >0$ for $z=re^{{\rm i}\theta}\in C$ with $\theta\in[-\pi+\eta,-\eta]$.
\\
On the other hand, we have $$\psi'(z)=(1-\lambda)^2\sqrt{\lambda}\frac{z^2-1}{(\sqrt{\lambda}z+1)^2(z+\sqrt{\lambda})^2}$$
so $\psi'$ vanishes only on  $1$ and $-1$, and for all $z$ in the unit circle such that  $z\neq\pm 1$, the gradient $\nabla(\Im \psi)(z)$ is non null and orthogonal to the tangent of the circle, hence \\
(e)  for $n$ large enough,  $\partial \Im \psi_n/\!\partial r>0$  on  $z=re^{{\rm i}\theta}\in A$ with $\theta\in[-\pi+\eta,-\eta]$. \\
At last, $\psi''(1)$ and $\psi''(-1)$ are non null, so \\
(f) for $n$ large enough,  $\psi_n'$ is one-to-one $\psi_n$ is at most two-to-one in a neighborhood of each of the disks $D=\{ z\; ;\; |z+1|\leq \alpha\}$ and $E=\{z\; ;\; |z-1|\leq\alpha\}$, 
\\
(g) for $n$ large enough,  $\psi_n'$ has no zero in both disks $\{z\; ;\; |z\pm 1|<\eta\}$.
\\
\\
For such a value of $n$, we actually have \\
(c') $\Im \psi_n(z) <0$ for $z=re^{{\rm i}\theta}\in B$ with $\theta\in]-\pi,0[$,\\
(d') $\Im \psi_n(z) >0$ for $z=re^{{\rm i}\theta}\in C$ with $\theta\in]-\pi,0[$.
\\
Indeed, if for instance $\Im \psi_n(z)$ were positive for some $z\in B$ in the lower half plane, then we would deduce the existence of   $z'\in B\cap\mathbb{C}^-$ for which $\Im \psi_n(z')=0$. But then $\psi_n(\overline{z'})=\overline{\psi_n(z')}=\psi_n(z')$, and this contradicts condition  (b). Let $u_n$ and $v_n$ be the unique zeros of  $\psi_n'$ in  $D$ and $E$, respectively. Observe that $u_n$ and $v_n$ must be real. Indeed, as $\psi_n'(\overline{u_n})=\overline{\psi_n'(u_n)}=0$ and $\overline{u_n}\in D$,   uniqueness implies $u_n=\overline{u_n}$. 
\\

Consider now the set $$F_{n,0}=\{ z\in\left(A\cup D\cup E \right)\; ;\;   z\notin\mathbb{R},\: \psi_n(z)\in\mathbb{R}\}\cup\{u_n,v_n\}.$$
 We can show that $F_{n,0}$ is an analytic curve that contains a curve  $\gamma_n$ joigning $u_n$ to $v_n$ in the half upper plane (the proof is the same as the one  page 219 in the article by D. Voiculescu and H. Bercovici (\cite{supercv})).  

Denote by $\Theta_n$ the part of the lower half plane below \\
$]-\infty,u_n]\cup \gamma_n\cup [v_n,+\infty[$. Let us show that $\psi_n$ is a one-to-one map from $\Theta_n$ onto the half upper plane. \\
For $R>100$, let $\gamma_{n,R}$ be the contour that runs (anticlockwise) through 
\begin{eqnarray*}&\gamma_n^*\:\cup\: [v_n,R]\:\cup\: \left\{z\in \mathbb{C}^-\; ;\; |z|=R\right\}\:\cup\: [-R,-\frac{1}{\sqrt{\lambda}}-\frac{1}{R}]&\\
&\cup \:\left\{z\in\mathbb{C}^-\; ;\; |z+\frac{1}{\sqrt{\lambda}}|=\frac{1}{R}\right\}\:\cup\: [-\frac{1}{\sqrt{\lambda}}+\frac{1}{R},u_n]. &\end{eqnarray*}
By the residues formula, it suffices to show that for every  $\omega$ in the upper half plane, the integral of $\frac{\psi_n'}{\psi_n-\omega}$ over the contour $\gamma_{n,R}$ tends to $2{\rm i}\pi$ when $R$ tends to infinity (this integral is well defined for  $R$  large enough, when $n$ and $\omega$ are fixed, because $\underset{\infty}{\lim}\psi_n=0$). 
\begin{itemize}
\item[\!\!\!\!\!\!(1)]  Computation of the limit of the integral on the complementary of the semi-circles in the curve $\gamma_{n,R}$:
\\
By our hypothesis on $f_n$, 
\begin{eqnarray*}\lim_{x\to-\infty}\psi_n(x)=0 &\textrm{ and } &\lim_{x\stackrel{<}{\to}-\frac{1}{\sqrt{\lambda}}}\psi_n(x)=+\infty,\\
\lim_{x\to+\infty}\psi_n(x)=0&\textrm{ and } &\lim_{x\stackrel{>}{\to}-\frac{1}{\sqrt{\lambda}}}\psi_n(x)=+\infty.
\end{eqnarray*}

The function $\psi_n$ is real on the complement of the semi-circles, so,  when $R$ tends to  $+\infty$,  the integral of $\frac{\psi_n'}{\psi_n-\omega}$ on this complement tends to $$\int_{-\infty}^{+\infty}\frac{\D y}{y-\omega}={\rm i}\pi.$$
\item[\!\!\!\!\!\!(2)] Computation of the limit of the integral on the large semi-circle:
\\
By hypothesis,  \[\lim_{R\to\infty}\left(R\sup_{|z|=R}\left | \frac{\psi_n'(z)}{\psi_n(z)-\omega}\right|\right)=0,\]
so the integral on the large semi-circle tends to  $0$ when $R$ tends to  infinity.
\item[\!\!\!\!\!\!(3)] Computation of the limit of the integral on the small semi-circle:
\\
By the same arguments as in  (1) and  (2), the integral of $\frac{\psi'}{\psi-\omega}$ on the complementary of the small semi-circle in the contour \begin{eqnarray*}&\left\{z\in\mathbb{C}^-\; ;\; |z|=1\right\}\cup [1,R]\cup \left\{z\in \mathbb{C}^-\; ;\; |z|=R\right\}\cup [-R,-\frac{1}{\sqrt{\lambda}}-\frac{1}{R}]&\\
&\cup \left\{z\in\mathbb{C}^-\; ;\; |z+\frac{1}{\sqrt{\lambda}}|=\frac{1}{R}\right\}\cup [-\frac{1}{\sqrt{\lambda}}+\frac{1}{R},-1] \end{eqnarray*} (running anticlockwise) tends to  ${\rm i}\pi$ when $R$ tends to  $+\infty$. Otherwise,  $\psi$ being a one-to-one map from  $\{z\in\mathbb{C}^-\; ;\; |z|>1\}$ to $\mathbb{C}^+$,   the integral of $\frac{\psi'}{\psi-\omega}$ on the  contour \begin{eqnarray*}&\left\{z\in\mathbb{C}^-\; ;\; |z|=1\right\}\cup [1,R]\cup \left\{z\in \mathbb{C}^-\; ;\; |z|=R\right\}\cup [-R,-\frac{1}{\sqrt{\lambda}}-\frac{1}{R}]&\\
&\cup \left\{z\in\mathbb{C}^-\; ;\; |z+\frac{1}{\sqrt{\lambda}}|=\frac{1}{R}\right\}\cup [-\frac{1}{\sqrt{\lambda}}+\frac{1}{R},-1]& \end{eqnarray*} (going anticlockwise) tends to  $2{\rm i}\pi$ when $R$ tends to  $+\infty$, so the integral of $\frac{\psi'}{\psi-\omega}$ on $\{z\in\mathbb{C}^-\; ;\; |z+1/\!\sqrt{\lambda}|=\frac{1}{R}\}$ tends to  $i\pi$ when $R$ tends to  $+\infty$. 
\\
But clearly, because of the boundness of $f_n$ and  $f_n'$ in a neighborhood of zero, we have $$\lim_{R\to+\infty}\int_{|z+\frac{1}{\sqrt{\lambda}}|,\Im z\leq 0}\left(\frac{\psi_n'(z)}{\psi_n(z)-\omega}-\frac{\psi'(z)}{\psi(z)-\omega}\right)\D z=0.$$
So the integral of $\frac{\psi_n'(z)}{\psi_n(z)-\omega}$ on the small semi-circle of  $\gamma_{n,R}$ tends to  $i\pi$ when $R$ tends to  $+\infty$.    
\end{itemize}
 
Thus, $\psi_n$ is a one-to-one map from $\Theta_n$ onto the upper half plane. Hence $K_n$ is a one-to-one map from $\{\frac{\sqrt{\lambda}}{1-\lambda}z+\frac{1}{1-\lambda}\; ;\;  z\in\Theta_n\}$, that we denote by $\Omega_n$,  onto $\mathbb{C}^+$. 

Let us show that the inverse function  $G_n$ of $K_n$ satisfies the hypothesis of the lemma.
\begin{itemize}
\item[\!\!\!\!\!\!(i)] Let us show that $zG_n(z)$ tends to $1$ when $|z|$ tends to  infinity. 
\\
$K_n$ is bounded in the complementary, in $\Omega_n$, of any open disk centered in zero with radius larger than $1+1/\!\sqrt{\lambda}$, so $G_n(z)$ tends to  zero when $|z|$ tends to  infinity.  But $zG_n(z)=K_n(G_n(z))G_n(z)$, so the limit of    $zG_n(z)$ when  $|z|$ tends to infinity is the limit of   $yK_n(y)$ when $y$ tends to  zero, that is $1$. 
\item[\!\!\!\!\!\!(ii)] The function that associates $K_n(z)-\frac{1}{z}$ to $z$  can be analytically extended to a neighborhood of zero in $\mathbb{C}$ because $f_n$ is analytic in a neighborhood of zero.
\end{itemize}
So the theorem is proved.
\end{proof}

The proof is very close to the proof given in \cite{supercv} of the failure of the Cram\'er theorem  for free random variables. In that article, the authors show also the superconvergence to the central limit theorem  for free random variables: in addition to the weak convergence, we know that for $n$ large enough,  the measures have an analytic density, and that their densities converge uniformly to the density of the semi-circle law. The superconvergence in the free Poisson limit theorem  has already been shown (using explicit computations) in \cite{hiai}.
\\
\\
{\bf Aknowledgments} 

The author would like to thank his advisor Philippe Biane for his comments on  draft of this paper.

Florent Benaych-Georges\\ DMA, \'Ecole Normale Sup\'erieure,\\ 45 rue d'Ulm, 75230 Paris Cedex 05, France\\  {\tt e-mail : benaych@dma.ens.fr}\\
 {\tt http://www.dma.ens.fr/$\sim$benaych}

\begin{thebibliography}{99}
\bibitem[A]{akhi} Akhiezer, N.I. \emph{The classical moment problem}, Moscou, 1961
\bibitem[BPB99]{appenice} Bercovici, H., Pata, V., with an appendix by Biane, P. \emph{Stable laws and domains of attraction in free
probability theory}, Annals of Mathematics, \textbf{149} (1999) 1023-1060
\bibitem[BV93]{defconv} Bercovici, H., Voiculescu, D. \emph{Free convolution of measures with unbounded supports}, Indiana Univ.
Math. J. \textbf{42} (1993) 733-773
\bibitem[BV95]{supercv} Bercovici, H., Voiculescu, D. \emph{Superconvergence to the central limit and failure of the Cram\'er theorem for free random variables}, Probability Theory and Related Fields \textbf{102} (1995) 215-222
\bibitem[HP00]{hiai} Hiai, F., Petz, D. \emph{The semicircle law, free random variables, and entropy}, Amer.
Math. Soc., Mathematical Surveys and Monographs Volume \textbf{77}, 2000
\bibitem[R]{rudin} Rudin, W. \emph{Real and complex analysis}, McGraw-Hill Book Co., New York, 1987.
\bibitem[S99]{spei2} Speicher, R. \emph{Notes of my lectures on Combinatorics of Free Probability (IHP, Paris, 1999)} available on \texttt{http://www.mast.queensu.ca/$\sim$speicher}
\bibitem[V86]{voic0}  Voiculescu, D.V. \emph{Addition of certain non-commuting random variables}, J. Funct. Anal. {\bf{66}} (1986) 323-346 
\bibitem[VDN92]{voic1} Voiculescu, D.V., Dykema, K., Nica, A. \emph{Free random variables}, CRM Monographs Series No.1, Amer. Math. Soc., Providence, RI, 1992 
\end{thebibliography}
\end{document}